\documentclass[a4]{article}
\usepackage{fullpage}
\usepackage{amsfonts}
\usepackage{amsmath}
\usepackage{graphicx}
\usepackage[mathscr]{eucal}

\title{On singularly perturbed linear cocyles over irrational rotations}
\author{Alexey V. Ivanov}
\date{}
\begin{document}
\renewcommand{\theequation}{\arabic{section}.\arabic{equation}}
\maketitle

\begin{abstract}
We study a linear cocycle over irrational rotation $\sigma_{\omega}(x) = x + \omega$  of a circle $\mathbb{T}^{1}$. It is supposed the cocycle is generated by a $C^{1}$-map 
$A_{\varepsilon}: \mathbb{T}^{1} \to SL(2, \mathbb{R})$ which depends on a small parameter $\varepsilon\ll 1$ and has the form of the Poincar\'e map corresponding to a singularly perturbed Schr\"odinger equation. Under assumption the eigenvalues of  $A_{\varepsilon}(x)$ to be of the form $\exp(\pm \lambda(x)/\varepsilon)$, where $\lambda(x)$ is a positive function, we examine the property of the cocycle to possess an exponential dichotomy (ED) with respect to the parameter $\varepsilon$. We show that in the limit $\varepsilon\to 0$ the cocycle "typically" exhibits ED only if it is exponentially close to a constant cocycle. In contrary, if the cocycle is not close to a constant one it does not posesses ED, whereas the Lyapunov exponent is "typically" large. 
\end{abstract}



\section{Introduction}

In the paper we study a skew-product map
\begin{equation}\label{eq1}
F_{A}: \mathbb{T}^{1}\times \mathbb{R}^{2} \to \mathbb{T}^{1}\times \mathbb{R}^{2}
\end{equation}
defined for any $(x, v) \in \mathbb{T}^{1}\times \mathbb{R}^{2}$ by 
$$
(x, v)\mapsto (\sigma_{\omega}(x), A(x)v),
$$
where $\sigma_{\omega}(x) = x + \omega$ is a rotation of a circle $\mathbb{T}^{1}$ with irrational rotation number $\omega$ and 
$$
A: \mathbb{T}^{1} \to SL(2, \mathbb{R})
$$ 
is a measurable function with respect to the Haar measure. 
The transformation $A$ generates a cocycle $M(x,n)$ by
\begin{eqnarray}
\nonumber
M(x, n)=A(\sigma_{\omega}^{n-1}(x))\ldots A(x), \,\, n>0;\\ 
\nonumber
M(x, n)=\left[A(\sigma_{\omega}^{-n}(x))\ldots A(\sigma_{\omega}^{-1}(x))\right]^{-1}, \,\, n<0;\\
\nonumber
\quad M(x, 0)=I.
\end{eqnarray}
Such discrete dynamical systems and their continuous counterparts are the subject of many papers appeared during the last three decades. Among important classes of examples of linear cocycles one needs to mention a discrete ergodic Schr\"odinger operator and quasiperiodic Hill's equation.
The main problems stated for this kind of systems include description of an asymptotic behaviour of their trajectories, establishing relations between dynamical characteristics and other (e.g. spectral) properties of a system, studying the genericity of systems with different behaviour in various classes of smoothness. Starting from the fundamental works \cite{DinSin}, \cite{JohnMos}, \cite{Her} there were elaborated many different methods  and techniques to investigate these problems \cite{Eli}, \cite{BourJit}, \cite{SorSpe}, \cite{BuFe}, \cite{AvBo} (see also \cite{DuaKle_book} and references therein). Nevertheless there are still many open questions. In particular, for a given family of cocycles the problem of effective description of a set of those parameter values which correspond to some specific property (e.g. hyperbolicity, reducibility, positiveness of the Lyapunov exponent) is far from its solution. The reason is a lack of constructive approaches. One of such approaches is based on the monodromization method developed in \cite{BuFe}, \cite{Fed}. This method was successfully applied to study the spectrum of the Harper operator and adiabatically perturbation of the 1-dimensional Schr\"odinger operator. Another constructive approach has been developed by different authors \cite{Jak}, \cite{BenCar}, \cite{LSY}, \cite{Laz} and is  based on inductive construction of the so-called "critical" set. We will discuss the details of this method in the end of the Section 2.

The paper is organized as follows. In Section 2 we specify the cocycle which will be studied and state the problem. Besides we remined some definitions and constructions. In particular, the idea of "critical set" method will be discussed. In Section 3 we describe the inductive procedure to construct the critical set and study its property. Finally, Section 4 is devoted to analysis of the Lyapunov exponent associated to the cocycle.

\section{Statement of the problem}
\setcounter{equation}{0}

In this paper we consider a one-parameter family of cocycles such that the corresponding  transformation $A$ has a special form. Particularly, given a $C^{2}$-function $f: \mathbb{T}^{2}\to \mathbb{R}$ of the two-torus, define sets
\begin{align*}
&\mathcal{D}_{+}=\{z\in \mathbb{T}^{2}: f(z)>0\},\\ 
&\mathcal{D}_{-}=\{z\in \mathbb{T}^{2}: f(z)<0\}
\end{align*}
and for any $x\in \mathbb{T}^{1}$ consider their intersections $S_{\pm}(x)=\mathcal{D}_{\pm}\cap I(x)$ with a segment $I(x) = \{(x + \omega s, s), s\in \mathbb{T}^{1}\}\subset \mathbb{T}^{2}$. We represent  $S_{\pm}(x)$ as
$$
S_{\pm}(x) = \bigcup\limits_{k=1}^{K(x)} \Delta_{k}^{\pm}(x),
$$
where $\Delta_{k}^{\pm}(x)$ are connected components ordered in a natural way with respect to increase of the parameter $s$ (see fig.1). 
\begin{figure}[h]
\center{\includegraphics[width=0.4\linewidth]{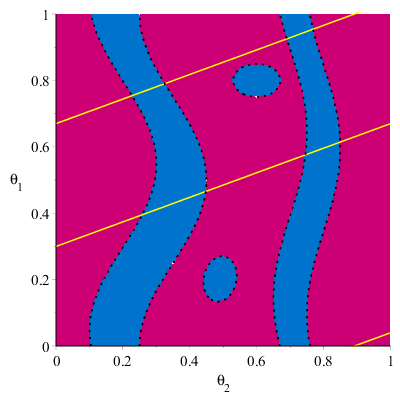}}
\caption{An example of two-torus with indicated domains $\mathcal{D}_{+}$ (red), $\mathcal{D}_{-}$ (blue) and segments $I(x)$ (yellow) drawn for two different values of $x=\theta_{1}$.}
\end{figure}

Let $g_{\pm}: \overline{\mathcal{D}_{\pm}} \to \mathbb{R}_{\pm}$ be $C^{2}$-functions.
Then the transformation $A$ is assumed to be of the form
\begin{equation}\label{eq2}
A(x) = \prod\limits_{k=0}^{K(x)}R(\varphi_{k}(x))\cdot Z(\lambda_{k}(x)), 
\end{equation}
where
\begin{equation*}
R(\varphi) = \left(
\begin{array}{cc}
\cos \varphi & \sin \varphi\\
-\sin \varphi & \cos \varphi
\end{array}
\right), 
\,
Z(\lambda) = \left(
\begin{array}{cc}
{\rm e}^{\lambda} & 0\\
0 & {\rm e}^{-\lambda}
\end{array}
\right)
\end{equation*}
and
\begin{equation}\label{eq3}
\varphi_{k}(x) = \varepsilon^{-1} \int\limits_{\Delta_{k}^{-}(x)}g_{-}{\rm d}l,\quad 
\lambda_{k}(x) = \varepsilon^{-1} \int\limits_{\Delta_{k}^{+}(x)}g_{+}{\rm d}l,\quad
\varepsilon\ll 1.
\end{equation}
We mention an evident property of the matrices $R$ and $Z$:
\begin{equation*}
R(\varphi_{1}+\varphi_{2}) = R(\varphi_{1})\cdot R(\varphi_{2}),\quad
Z(\lambda_{1}+\lambda_{2}) = Z(\lambda_{1})\cdot Z(\lambda_{2}),\quad
\forall\; \varphi_{1}, \varphi_{2}, \lambda_{1}, \lambda_{2}.
\end{equation*}

To emphasize the dependence of the transformation $A$ on the parameter $\varepsilon$ we will denote it as $A_{\varepsilon}(x)$.  Such linear cocycle  (\ref{eq1}) appears as a model in the problem of stability for the singularly perturbed Schr\"odinger equation with quasi-periodic potential $f$. Indeed, consider an equation 
\begin{equation}\label{eq4}
\varepsilon^{2}y'' = f(t)y
\end{equation}
 and assume the potential $f(t) = F(\theta_{0} + \omega t, t)$ where $F: \mathbb{T}^{2}\to \mathbb{R}$, i.e $f$ is a quasi-periodic function.
In the extended phase space this equation can be written as
\begin{equation*}
\varepsilon^{2}y'' = F(\theta_{1}, \theta_{2})y,\quad 
\theta'_{1} = \omega, 
\quad \theta'_{2} = 1.
\end{equation*}
One may consider the Poincar\'e map associated with this system: 
\begin{eqnarray}
\nonumber
\Phi: \left(
\begin{array}{cc}
y(0) \\
y'(0)\\
\theta_{1}(0)\end{array}
\right) 
\rightarrow 
\left(
\begin{array}{cc}
y(1) \\
y'(1)\\
\theta_{1}(1)\end{array}
\right) \,.
\end{eqnarray}
which is equivalent to
\begin{eqnarray}
\nonumber
\theta_{1}(1) = \theta_{1}(0) + \omega,\\
\label{eq5}
\left(
\begin{array}{cc}
y(1) \\
y'(1)\end{array}
\right)  
= A_{\varepsilon}(\theta_{1}(0)) 
\left(
\begin{array}{cc}
y(0) \\
y'(0)\end{array}
\right)
\end{eqnarray}
with some $A_{\varepsilon}(\theta)\in SL(2,\mathbb{R})$.
Thus, we arrive at a system of the form (\ref{eq1}). Using WKB-method one may show that the matrix $A_{\varepsilon}$ in (\ref{eq5}) can be written in the form (\ref{eq2}), where the functions $g_{\pm}$ are restrictions of $\vert f\vert^{1/2}$ on $\mathcal{D}_{\pm}$. If the potential $f$ is periodic then the stability of the origin for the system (\ref{eq4}) does not depend on the value $\theta_{1}(0)$ and the following theorem holds \cite{Iva17}
\newtheorem{thms}{Theorem}
\begin{thms}
There exists $\varepsilon_{0}>0$ and a subset 
$\mathcal{E}_{h} \subset (0, \varepsilon_{0})$ such that 
\newline 
\newline
1. for any $\varepsilon_{1}<\varepsilon_{0}$ the Lebesgue measure $leb\left((0,\varepsilon_{1})\setminus \mathcal{E}_{h}\right) = O\left({\rm e}^{-c/ \varepsilon_{1}}\right)$ with some positive constant $c$;
\newline 
\newline
2. for any $\varepsilon\in \mathcal{E}_{h}$ the origin is an unstable equilibrium of the system (\ref{eq3}).
\end{thms}
Note also that instability of the equilibrium is related to the exponential dichotomy possessed by the system (\ref{eq4}) for $\varepsilon\in \mathcal{E}_{h}$ \cite{Iva17}.

Before stating the results we remind the definition of the exponential dichotomy and some other concepts in cocycle settings.
\newtheorem{defs}{Definition}
\begin{defs}
A cocycle $M$ is said to have an exponential dichotomy (ED) if there are positive constants $K, \gamma$ and a projector-valued function $P(x)$ continuously dependent on $x\in \mathbb{T}^{1}$ such that
\begin{align}
\nonumber
\left\Vert M(x, m)P(x)M^{-1}(x, n)\right\Vert \le K{\rm e}^{-\gamma(m-n)}, \quad
m\ge n;\\
\nonumber
\left\Vert M(x, m)\bigl(I-P(x)\bigr)M^{-1}(x, n)\right\Vert \le K{\rm e}^{-\gamma(m-n)}, \quad
m\le n.
\end{align}
\end{defs}
One has to note that property to possess ED is very robust under perturbation as follows from \cite{Cop}, \cite{Palm84}, \cite{Palm87}.

A closely related to the concept of ED is the concept of uniform hyperbolicity.
\begin{defs}
A cocycle $M$ is said to be uniformly hyperbolic (UH) if there exist continuous maps $E^{u,s}: \mathbb{T}^{1} \to Gr(2,1)$ and positive constants $C, \Lambda$  such that the subspaces 
$E^{u,s}(x)$ are invariant with respect to the map (\ref{eq1})  
(i.e. $E^{u,s}(\sigma_{\omega}(x))=A(x)E^{u,s}(x)$) and $\forall x\in \mathbb{T}^{1}$, $n\ge 0$
\begin{align}
\left\Vert M(x, -n)\vert_{E^{u}(x)}\right\Vert \le C{\rm e}^{-\Lambda n}, \notag\\
\left\Vert M(x, n)\vert_{E^{s}(x)}\right\Vert \le C{\rm e}^{-\Lambda n}.\notag
\end{align}
\end{defs}
It can be shown (see e.g. \cite{AvBo}, \cite{Fab}) that a cocycle $M$ possesses ED if and only if it is UH. Moreover \cite{AvBo}, these properties are equivalent to existence of  positive constants $C, \Lambda$ such that $\forall x\in \mathbb{T}^{1}$ and $n\ge 0$
$$\Vert M(x,n)\Vert \ge C {\rm e}^{\Lambda n}.$$
Note that rotation $\sigma_{\omega}$ is ergodic with respect to Haar measure since $\omega$ is assumed to be irrational. This leads by Oseledets theorem to existence of invariant subspaces $E^{u,s}(x)$ for a.e. $x\in \mathbb{T}^{1}$. However, in general, the maps $E^{u,s}$ are only measurable.
Besides, by Kingman's subadditive ergodic theorem for a.e. $x\in \mathbb{T}^{1}$ there exists the Lyapunov exponent 
$$\Lambda(x) = \lim\limits_{n\to \infty}\frac{1}{n}\log\Vert  M(x,n)\Vert$$
and due to ergodicity of $\sigma_{\omega}$
$$\Lambda(x) = \Lambda_{0}\quad a.e.,$$
where $\Lambda_{0}$ is the integrated Lyapunov exponent 
$$\Lambda_{0} = \int \Lambda(x){\rm d}x = \lim\limits_{n\to \infty}\int \frac{1}{n}\log\Vert  M(x,n)\Vert{\rm d}x.$$

We also remind the definition of reducibility.
\begin{defs}
A cocycle $M$ is said to be reducible if it is conjugated to a constant cocycle, i.e. there exists a continuous map $B: \mathbb{T}^{1}\to {\rm SL}(2, \mathbb{R})$ and 
a constant $C\in {\rm SL}(2, \mathbb{R})$ such that
$$
B(\sigma_{\alpha}(x)) A(x) B^{-1}(x) = C, \quad \forall x\in \mathbb{T}^{1}.
$$
\end{defs}
That is a cocycle is reducible if it is conjugated to a constant cocycle by a suitable coordinate change.

All these properties are related by the following theorem (see e.g. \cite{Avi}, \cite{AvBo})
\begin{thms}
Let $M$ be a cocycle defined by the skew-product map (\ref{eq1}). Then the following statements are equivalent
\newline
1. $M$ possesses ED; 
\newline
2. $M$ is UH;
\newline
3. $M$ is reducible and $\Lambda_{0}>0$.  
\end{thms}

In the present work, using the approach developed in \cite{Jak}, \cite{BenCar}, \cite{LSY},\cite{Laz}, we analyze those values of the parameter $\varepsilon$ for which the cocycle associated to the map (\ref{eq1}) is reducible and has positive Lyapunov exponent. Although papers \cite{Jak}, \cite{BenCar}, \cite{LSY},\cite{Laz} are devoted to different objects, they use a common framework. The idea of the method suggested in \cite{Jak} can be described in the following way. Consider a one-parameter family of skew-product maps $F_{A,\varepsilon} = (f, A_{\varepsilon})$ defined (similarly to (\ref{eq1})) on a vector bundle $V$ over a base $B$, then properties  of the fiber transformation $A_{\varepsilon}$ depend (in general) on a point of the base. Choose those points of the base which correspond to the violation of some specific property (e.g. hyperbolicity) of $A_{\varepsilon}$. These points constitute an initial approximation $\mathcal{C}_{0}$ of the critical set. Taking a small neighborhood of $\mathcal{C}_{0}$ one may study the dynamics of this set under the map $f$. Interactions between different parts of $\mathcal{C}_{0}$ give rise to the first approximation $\mathcal{C}_{1}$ of the critical set. Detuning the parameter $\varepsilon$ one may put these interactions in "general" position to move the induction forward and to control properties of $\mathcal{C}_{1}$. Continueing the induction one may hope to construct the critical set as a limit of $\mathcal{C}_{n}$ and investigate a set $X_{h}$ of those points of the base which do not approach this critical set too closely for sufficiently large period of time. It might happen that the set $X_{h}$ is sufficiently large. Then, using properties of the map $f$ (e.g. ergodicity), one may extract an additional information on the whole system. We follow this approach in the present paper. The results can be formulated in two theorems:

\begin{thms} 
Assume $\min\limits_{k, x} \lambda_{k}(x) = \lambda_{0}/\varepsilon>0$ and for any $k$ $\varphi_{k}(x) = \phi_{k}/\varepsilon$, where $\varphi_{k} = const$. Then there exist sufficiently small $\varepsilon_{0}>0$, positive constants $C_{\Lambda}<1$, $C_{0}$  and a subset $\mathcal{E}_{h} \subset ( 0, \varepsilon_{0})$ such that 
\newline 
\newline
1. for any $\varepsilon_{1}<\varepsilon_{0}$ the Lebesgue measure 
$leb\left((0, \varepsilon_{1})\setminus \mathcal{E}_{h}\right) = 
O\left({\rm e}^{-C_{0}/ \varepsilon_{1}}\right)$;
\newline 
\newline
2. for any $\varepsilon\in\mathcal{E}_{h}$ the Lyapunov exponent $\Lambda_{0} > {\rm e}^{C_{\Lambda}\lambda_{0}/\varepsilon}$;
\newline 
\newline
3. for any $\varepsilon\in \mathcal{E}_{h}$ the cocycle $M$ possesses ED. 
\end{thms}

This theorem is an analog of Theorem 1. Indeed, if $\omega$ is rational and, hence, the function $F$ associated to (\ref{eq4}) is periodic, one may always assume without loss of generality that $F$ does not depend on $\theta_{1}$. It means the functions $\varphi_{k}$, which appear in representation of the matrix $A$ in (\ref{eq5}) in the form (\ref{eq2}), are also independent of $\theta_{1}$. Note that in this case $\lambda_{k}$ become constant too. It will be shown later that dependence of $\lambda_{k}$ on $x$ is not essential whereas inconstancy of $\varphi_{k}$ is crucial.

The second theorem describes the opposite case when the functions $\varphi_{k}$ are not constant. Let $p_{n}/q_{n}$ denotes the rational approximation of order $n$ to $\omega$ in its continued fraction expansion.
\begin{defs}
We will say that $\omega$ satisfies the Brjuno's condition with a constant $C_{B}$  if
$$\sum\limits_{n=1}^{\infty}\frac{\log (2 q_{n+1})}{q_{n}} = C_{B} < \infty. $$
\end{defs}
This condition was introduced by A. D. Brjuno in \cite{Brj} and also used in \cite{LSY}.  In this paper we consider additional condition:
\begin{defs}
Assume $\omega$ satisfies the Brjuno's condition with a constant $C_{B}$. We will say that $\omega$ satisfies a condition $(A)$ if there exist a subsequence $\{q_{n_{j}}\}_{j=1}^{\infty}$ and positive constants $C_{\omega}$,  $C_{\varepsilon}$, $C_{\delta}$, $\gamma$ such that $C_{\delta}<1$
\begin{equation}\label{eq6}
q_{n_{j}+1} > C_{\omega} q_{n_{j}}^{1+\gamma},\;\; \forall\; j
\end{equation}
and for all $k\in \mathbb{N}$ there exists an index $J_{k}$
\begin{equation}\label{eq7}
\frac{1}{q_{n_{J_{k}}}}\left(\log q_{n_{J_{k}}} + \frac{\log C_{\delta}^{-1}}{1+\gamma}\right)<\frac{\log q_{n_{J_{k+1}}}}{q_{n_{J_{k}}}} < \\
C_{\varepsilon} - \frac{C_{B}}{1+\gamma}\left[1-C_{\delta}-\frac{1}{q_{n_{J_{k}}}}\right] - 
\frac{\log C_{\omega}}{(1+\gamma)q_{n_{J_{k}}}}.
\end{equation}
\end{defs}
The inequality (\ref{eq6}) says that for the sequence of denominators $\{q_{n}\}$ the ratio $q_{n+1}/q_{n}$ becomes sufficiently large infinitely many times, while the inequality (\ref{eq7}) guarantees that such events of growth of the denominators occur in a regular way.
\begin{thms} 
Assume $\min\limits_{k, x} \lambda_{k}(x) = \lambda_{0}>0$ and for any $k$  the function $\varphi_{k}$ does not have degenerate critical points. If
$\omega$ satisfies the Brjuno's condition with constant $C_{B}$ and the condition $(A)$ 
then there exist sufficiently small $\varepsilon_{0}>0$, positive constants $C_{\Lambda}<1$, $C_{0}$  and a subset $\mathcal{E}_{h} \subset ( 0, \varepsilon_{0})$ such that 
\newline 
\newline
1. for any $\varepsilon_{1}<\varepsilon_{0}$ the Lebesgue measure 
$leb\left((0, \varepsilon_{1})\setminus \mathcal{E}_{h}\right) = 
O\left({\rm e}^{-C_{0}/ \varepsilon_{1}}\right)$;
\newline 
\newline
2. for any $\varepsilon\in\mathcal{E}_{h}$ the Lyapunov exponent $\Lambda_{0} > {\rm e}^{C_{\Lambda}\lambda_{0}/\varepsilon}$;
\newline 
\newline
3. for any $\varepsilon\in (0, \varepsilon_{0})$ the cocycle $M$ does NOT possess ED. 
\end{thms}

{\bf Remark} One has to remark that Theorems 3, 4 are formulated in terms of the functions $\varphi_{k}$, $\lambda_{k}$ only. It is essential their dependence on the parameter $\varepsilon$, but the representation via integrals (see (\ref{eq3})) is not essential. Thus, instead of this representation we ma consider the following one $\varphi_{k} = \varepsilon^{-1}\hat\varphi_{k}$, $\lambda_{k} = \varepsilon^{-1}\hat\lambda_{k}$, with $\hat\varphi_{k}$, $\hat\lambda_{k}$ to be real-valued $C^{2}$-functions defined on a segment $I_{k}\subset \mathbb{T}^{1}$. However, we prefer the definition (\ref{eq2}), since it reveals the nature of the problem.

\section{Critical set}
\setcounter{equation}{0}

In this section we describe an inductive procedure for the construction of a critical set consisting of those points $x\in \mathbb{T}^{1}$ where the cocycle $M$ exhibits "bad" hyperbolic behaviour.

Given two sequences of positive numbers $\{\varphi_{k}\}, \{\lambda_{k}\}_{ k=1}^{\infty}$ we introduce $A_{k} = R(\varphi_{k})\cdot Z(\lambda_{k})$ and $A^{k} = \prod\limits_{j=1}^{k} A_{j}$. Then the transformation (\ref{eq2}) takes the form 
$A_{\varepsilon}(x) =  A^{K(x)}(x)$, where the sequences $\{\varphi_{k}\}, \{\lambda_{k}\}$ are defined by (\ref{eq3}). Note that although $\Vert A_{k}\Vert = \lambda_{k}$,  it may happen that $\Vert A_{k}\cdot A_{j}\Vert \ll \lambda_{k}\lambda_{j}$. Nevertheless, one may prove the following lemma.

\newtheorem{lemmas}{Lemma}
\begin{lemmas}
Let $\{\varphi_{k}\}, \{\lambda_{k}\}_{ k=1}^{\infty}$ be two sequences of positive numbers such that  for any $k$
$$
\vert \cos\varphi_{k}\vert \ge \delta,\quad \lambda_{k}\ge \lambda_{0}>0,\quad
\delta {\rm e}^{\lambda_{0}}\gg 1.
$$
Then 
$$\Vert A^{k}\Vert \ge \left(C_{A}\delta {\rm e}^{\lambda_{0}}\right)^{k},\, \forall\, k,$$
where the constant $C_{A} = 1 + O\left(\delta {\rm e}^{\lambda_{0}}\right)^{-2}$.
\newline
Moreover, $A^{k}$ can be represented as
\begin{equation*}
A^{k} = R(\theta_{k}+\chi_{k})Z(\mu_{k})R(-\chi_{k})
\end{equation*}
with 
$$
\mu_{k}\ge k(\lambda_{0} + {\rm ln}\delta - {\rm ln}C_{A}), \quad
\chi_{k} = O\left(\delta^{-1}{\rm e}^{-2\lambda_{0}}\right).
$$
\end{lemmas}

PROOF: Let $\alpha_{1}, \alpha_{2}, \varphi$ be real numbers such that $\alpha_{j} > 1, j=1,2$. Then, using the polar decomposition one may represent
\begin{equation*}
Z(\log\alpha_{2}) R(\varphi), Z(\log\alpha_{1}) = R(\theta+\chi)Z(\mu)R(-\chi),
\end{equation*}
where
\begin{multline}\label{eq8}
\cos\theta = \frac{\cot\varphi}{\left(\cot^{2}\varphi + 
\left(
\frac{\alpha_{1}^{-2} + \alpha_{2}^{-2}}{1 + \alpha_{1}^{-2}\alpha_{2}^{-2}}
\right)^{2}\right)^{1/2}},\quad
\sin\theta = \frac{
\frac{\alpha_{1}^{-2} + \alpha_{2}^{-2}}{1 + \alpha_{1}^{-2}\alpha_{2}^{-2}}
}{\left(\cot^{2}\varphi + 
\biggl( 
\frac{\alpha_{1}^{-2} + \alpha_{2}^{-2}}{1 + \alpha_{1}^{-2}\alpha_{2}^{-2}} \biggr)^{2} \right)^{1/2}},\\
{\rm e}^{\mu} = \frac{a}{2}\left(1+ c + \left(\left(1+ c\right)^{2}-4 a^{-2}\right)^{1/2}\right),\quad
a=\alpha_{1}\alpha_{2}\cos\varphi \cos\theta \left(1 + \alpha_{2}^{-2} \frac{\alpha_{1}^{-2} + \alpha_{2}^{-2}}{1 + \alpha_{1}^{-2}\alpha_{2}^{-2}} \tan^{2}\varphi\right),\\
b=\alpha_{1}^{-2}\tan\varphi\frac{1-\alpha_{2}^{-2}}
{1+\alpha_{1}^{-2}\alpha_{2}^{-2} + 
\alpha_{2}^{-2}(\alpha_{1}^{-2}+\alpha_{2}^{-2})\tan^{2}\varphi},\quad
c=b^{2}+a^{-2},\\
\cos(2\chi) = \frac{1-c}{(1-c)^{2}+4 b^{2}},\quad
\sin(2\chi) = -\frac{2 b}{(1-c)^{2}+4 b^{2}}.
\end{multline}
Taking this into account we obtain clearly for $k=1$ that $A^{1}=A_{1} =R(\theta_{1}+\chi_{1})Z(\mu_{1})R(-\chi_{1})$ with
$$
\theta_{1} = \varphi_{1},\quad
\chi_{1} = 0,\quad
\mu_{1}\ge \lambda_{0}.
$$ 
Applying consecutively  (\ref{eq8}) one gets for $k=n>1$
\begin{multline}\label{eq9}
\vert\theta_{n} - \varphi_{n}\vert \le 2\delta^{-1}{\rm e}^{-2\lambda_{0}}, \\
\vert\chi_{n} - \chi_{n-1}\vert \le \delta^{-2n+3}{\rm e}^{-2(n-1)\lambda_{0}}
\left(1-2\delta^{-2}{\rm e}^{-4\lambda_{0}}\right)^{-2(n-2)}, \\
\mu_{n} - \mu_{n-1} \ge \lambda_{0} + \log\delta - 2\delta^{-2}{\rm e}^{-4\lambda_{0}}. 
\end{multline}
Summation over $n$ finishes the proof. $\square$

This lemma suggests to define the initial approximation of the critical set as
$$
\mathcal{C}_{0}(\varepsilon) = \left\{x\in \mathbb{T}^{1}: \; \exists\,  (k, i): \; \varphi_{k}(x) = \frac{\pi}{2}(1 + 2i)\varepsilon\right\}.$$ 
Here $k\ge 1, i\ge 0$. This set consists of finite number $N(\varepsilon)$ of points, which we denote by $c_{j}^{(0)}(\varepsilon)$. Thus,
$$
\mathcal{C}_{0}(\varepsilon) = \bigcup\limits_{j=1}^{N(\varepsilon)}\{c_{j}^{(0)}(\varepsilon)\}.
$$ 
Let $\rho_{j,j'}(\varepsilon) = dist(c_{j}^{(0)}(\varepsilon), c_{j'}^{(0)}(\varepsilon))$, where $dist(x_{1}, x_{2})$ stands for the distance between $x_{1}$ and $x_{2}$. In the rest of the paper we will assume the functions $\varphi_{k}$ to be such that 
\begin{equation}\label{eq10}
\left\vert\frac{d\rho_{j,j'}}{d\varepsilon}\right\vert\ge C_{\rho}>0.
\end{equation}
We consider  
$$
\mathcal{E}_{p,0} = \{\varepsilon > 0: N(\varepsilon) = p\} = \bigcup (\varepsilon_{p, 1}, \varepsilon_{p, 0}]
$$
and for $\varepsilon\in \mathcal{E}_{p,0}$ introduce a set which will be called (following \cite{Laz}) a layer
\begin{equation*}
L\mathcal{C}_{0} = 
\bigcup\limits_{j=1}^{N(\varepsilon)} (c_{j}^{(0)}-\delta_{0}, c_{j}^{(0)}+\delta_{0}) = \bigcup\limits_{j=1}^{N(\varepsilon)} I_{j,0}, \quad 
\delta_{0} = {\rm e}^{-\lambda_{0}\kappa_{0}/\varepsilon_{p,0}},
\end{equation*}
where $\kappa_{0}$ will be chosen later.

Fix a point $x\in \mathbb{T}^{1}$ and denote its finite trajectory under the map $\sigma_{\omega}$ by $x_{k}=\sigma_{\omega}^{k}(x), k=0,\ldots, m$. Note that if this trajectory does not fall into $L\mathcal{C}_{0}$ one may apply Lemma 1 to the cocycle $M(x, k)$ and obtain the following estimate
\begin{equation*}
M(x,k) = R(\theta_{k}(x)+\chi_{k}(x))Z(\mu_{k}(x))R(-\chi_{k}(x))
\end{equation*}
with
\begin{equation}\label{eq11}
\mu_{k}(x) \ge \frac{k}{\varepsilon}\Bigl(\lambda_{0} - \kappa_{0} - \log C_{A}\Bigr),\quad
\chi_{k}(x) = O\left(\varepsilon {\rm e}^{-(2-\kappa_{0})\lambda_{0}/\varepsilon}\right),\quad
\forall\; 0\le k\le m.
\end{equation}
 However, due to irrationality of $\omega$ for any $x\in \mathbb{T}^{1}$ there exists a time $\tau(x)$ such that $\sigma_{\omega}^{\tau(x)}(x)\in L\mathcal{C}_{0}$. Consider the dynamics of the layer $L\mathcal{C}_{0}$ itself. 
\begin{defs}
Let $\tau_{j,j'}^{(0)}$ be the minimum of $k\ge 1$ such that
$$
\sigma_{\omega}^{k}(I_{j,0})\cap I_{j',0} \neq \emptyset.
$$
We call such event a collision of the zeroth order and $\tau_{j,j'}^{(0)}$ the time of collision.  
A collision is called primary if $j=j'$ and secondary if $j\neq j'$. 
\end{defs}
One may note that the only possibility to avoid collisions occurs in a trivial case when the set $\mathcal{C}_{0}(\varepsilon)$ is empty. If the functions $\varphi_{k}$ are not constant the number of critical points corresponding to some particular $k_{*}$ can be estimated from below by the integer part of $\frac{1}{\pi\varepsilon}Var(\varphi_{k_{*}})$, where $Var$ denotes the variance of a function. Hence, $\mathcal{C}_{0}(\varepsilon)\neq \emptyset$ for sufficiently small $\varepsilon$. On the other hand, if all $\varphi_{k}$ are constant one may exclude a countable set of values of the parameter $\varepsilon$ to get $\mathcal{C}_{0}(\varepsilon)= \emptyset$. We exclude a larger set to apply Lemma 1. Namely, let
$\varphi_{k}(x)=\phi_{k}/\varepsilon$, $\phi_{k}=const$. Define $\varepsilon_{k,j} = \pi\phi_{k}(1/2 +j)^{-1}, 
\delta_{k,j} = \varepsilon_{k,j}\phi_{k}^{-1}{\rm e}^{-\lambda_{0}/2\varepsilon_{k,j}}$ and
$$
\mathcal{E}_{e} = \bigcup\limits_{k,j}(\varepsilon_{k,j}-\delta_{k,j}, \varepsilon_{k,j}+\delta_{k,j}).
$$
Then for any $\varepsilon_{0}>0$ the Lebesgue measure
$$
leb(\mathcal{E}_{e}\cap (0, \varepsilon_{0})) = 
O\left({\rm e}^{-\lambda_{0}/2\varepsilon_{0}}\right)
$$
and for any $\varepsilon\in \mathcal{E}_{h}=(0, \varepsilon_{0})\setminus \mathcal{E}_{e}$ one has $\vert\cos\varphi_{k}(x)\vert \ge {\rm e}^{-\lambda_{0}/2\varepsilon_{0}}$. Thus, applying Lemma 1 to the cocycle $M$, corresponding to $\varepsilon\in \mathcal{E}_{h}$, we obtain Theorem 3.

Note that if $\varphi_{k}$ are not constant, but thier variances are sufiiciently small then the set $\mathcal{C}_{0}(\varepsilon)$ might be empty for intermediate values of $\varepsilon$. Particularly, one gets the following corollary of Theorem 3
\newtheorem{crls}{Corollary}
\begin{crls}
Let functions $\lambda_{k}$, constants $\phi_{k}$, $C_{\Lambda}$, $C_{0}$ and the parameter $\varepsilon_{0}$ satisfy Theorem 3.  Then for a cocycle defined by the transformation (\ref{eq2}) associated to the set of functions $\lambda_{k}(x)$ and 
$\varphi_{k}(x) = (\phi_{k} + \psi_{k}(x))/\varepsilon$ such that
$$
Var(\psi_{k})\ll {\rm e}^{-\lambda_{0}/2\varepsilon_{0}},\;\; \forall k
$$
there exist $\varepsilon_{1}\ll \varepsilon_{0}$ and a subset $\mathcal{E}_{h} \subset 
(\varepsilon_{1}, \varepsilon_{0})$ such that 
\newline 
\newline
1. the Lebesgue measure 
$leb\left((\varepsilon_{1}, \varepsilon_{0})\setminus \mathcal{E}_{h}\right) = 
O\left({\rm e}^{-C_{0}/ \varepsilon_{0}}\right)$;
\newline 
\newline
2. for any $\varepsilon\in\mathcal{E}_{h}$ the Lyapunov exponent $\Lambda_{0} > {\rm e}^{C_{\Lambda}\lambda_{0}/\varepsilon}$;
\newline 
\newline
3. for any $\varepsilon\in \mathcal{E}_{h}$ the cocycle $M$ possesses ED.
\end{crls}

Now let us assume that $\varphi_{k}$ are not constant and $\varepsilon$ is sufficiently small such that $\mathcal{C}_{0}(\varepsilon)\neq \emptyset$ . In this case we cannot avoid the collisions and our goal is to make the collisions as rare as possible. It is evident that primary and secondary collisions behave in a different manner with respect to varying the parameter $\varepsilon$. Indeed, the times of primary collisions depend on $\varepsilon$ via $\delta_{0}$ only and are characterized essentially by the rotation number $\omega$.  Whereas under assumption (\ref{eq10}) one may hope to detune $\varepsilon$ such that the times of secondary collisions would be greater than the times of primary ones. We also remark that for all $j$ the times $\tau_{j,j}^{(0)}$ are equal and we denote them by $\tau_{0}$. Define 
$$
\mathcal{E}'_{p,0}=\{\varepsilon\in \mathcal{E}_{p,0}: \exists\; (j,j'):\; \tau_{j,j'}^{(0)}<\tau_{0}\}
$$
and estimate the Lebesgue measure of this set. If $\{q_{n}\}$ denote the denominators of the best rational approximations for $\omega$ then $q_{1}, q_{2},\dots$ are exactly the times when $\sigma_{\omega}^{m}(x)$ approximate $x$ better than ever before. Besides, the following inequality holds
\begin{equation}\label{eq12}
\frac{1}{2q_{n+1}} < dist(\sigma_{\omega}^{q_{n}}(x),x) < \frac{1}{q_{n+1}},\quad
\forall\; n.
\end{equation} 
In the rest of the paper we will suppose that $\omega$ satisfies the Brjuno's condition and the condition $(A)$.  Under this assumption we choose a sufficiently large $J_{0}$ and $\kappa_{0}$ such that
\begin{equation}\label{eq13}
\frac{1}{q_{k_{J_{0}}+1}} < \delta_{0} = {\rm e}^{-\lambda_{0}\kappa_{0}/\varepsilon} = \frac{1}{C_{\omega}q_{k_{J_{0}}}^{1+\gamma}}.
\end{equation}
This choice and (\ref{eq12}) imply
\begin{eqnarray}\label{eq14}
\nonumber
dist\left(\sigma_{\omega}^{m}(x), x\right) > \frac{1}{2 q_{n_{J_{0}}}}\gg 
\frac{1}{C_{\omega}q_{n_{J_{0}}}^{1+\gamma}}=\delta_{0},\; 0\le m\le q_{n_{J_{0}}}-1,\\
dist\left(\sigma_{\omega}^{q_{n_{J_{0}}}}(x), x\right) < \frac{1}{q_{n_{J_{0}}+1}} < 
\frac{1}{C_{\omega}q_{n_{J_{0}}}^{1+\gamma}}=\delta_{0}.
\end{eqnarray}
Hence, one may conclude that $\tau_{0} = q_{n_{J_{0}}}$. To estimate $leb(\mathcal{E}'_{p,0})$ note that 
\begin{equation}\label{eq15}
N(\varepsilon) = O\bigl(\varepsilon^{-1}\bigr),\quad 
\varepsilon_{p,0}-\varepsilon_{p,1} = O\bigl(\varepsilon_{p,0}^{2}\bigr)
\end{equation}
and due to (\ref{eq10}) $\Delta\rho = \max \rho_{j,j'}(\varepsilon) - \min \rho_{j,j'}(\varepsilon)$ (where $\max$ and $\min$ are taken over $\varepsilon\in (\varepsilon_{p,1}, \varepsilon_{p,0}]$ and all $(j,j')$) is
\begin{equation}\label{eq16}
\Delta\rho = O(\varepsilon_{p,0}^{2}).
\end{equation}
Besides the number of collisions is estimated by $\left(\begin{array}{c}2\\p\end{array}\right)$ and we get
\begin{equation}\label{eq17}
leb(\mathcal{E}'_{p,0})\le C_{\rho}\Delta\rho \left(\begin{array}{c}2\\p\end{array}\right)
\tau_{0}\delta_{0}.
\end{equation}
Substituting (\ref{eq13}), (\ref{eq15}), (\ref{eq16}) into (\ref{eq17}) we arrive at the following lemma.
\begin{lemmas}
The Lebesgue measure of $\mathcal{E}'_{p,0}$ is estimated as
$$ 
leb(\mathcal{E}'_{p,0}) = 
O\left({\rm e}^{-\lambda_{0}\kappa_{0}\gamma/(1+\gamma)\varepsilon_{p,0}}\right).
$$
\end{lemmas}
We eliminate this "bad" set and introduce 
$\hat{\mathcal{E}}_{p,1} = \mathcal{E}_{p,0}\setminus \mathcal{E}'_{p,0}$.

For $\varepsilon\in \hat{\mathcal{E}}_{p,1}$ and $x\in I_{j,0}$ we consider $M(x, \tau_{0})$. Applying Lemma 1 to 
$M(A(x,\varepsilon),\tau_{0}-1)$ and using (\ref{eq11}) we represent it as
\begin{equation}\label{eq18}
M(A(x,\varepsilon),\tau_{0}-1) = 
R(\theta_{j,0}(x)+\chi_{j,0}(x))Z(\mu_{j,0}(x))R(-\chi_{j,0}(x)),
\end{equation}
with 
\begin{equation*}
\chi_{j,0}(x) = O(\varepsilon{\rm e}^{-(2-\kappa_{0})\lambda_{0}/\varepsilon_{p,0}}),\quad
\mu_{j,0}(x)\ge \frac{\tau_{0}-1}{\varepsilon}(\lambda_{0}-\kappa_{0}-\log C_{A}).
\end{equation*} 
We substitute (\ref{eq18}) into the expression for $M(x, \tau_{0})$ to get
\begin{equation*}
M(x,\tau_{0}) = 
R(\theta_{j,0}(x)+\chi_{j,0}(x))Z(\mu_{j,0}(x))R(\varphi_{k}(x)-\chi_{j,0}(x))
Z(\lambda_{k}(x)),
\end{equation*}
with $k$ defined by
\begin{equation}\label{eq19}
\cos(\varphi_{k}(c_{j}^{(0)})) = 0.
\end{equation}

Then we consider an equation
\begin{equation}\label{eq20}
\cos\left(\varphi_{k}(x) - \chi_{j,0}(x)\right) = 0,\quad x\in I_{j,0},
\end{equation}
where $k$ satisfies (\ref{eq19}).  One notes that if a point $x$ is close to solution of (\ref{eq20}), $M(x, \tau_{0})$ loses the hyperbolicity. 

Expanding (\ref{eq20}) around $x=c_{j}^{(0)}$ we obtain an equation for 
$\Delta_{j} = x- c_{j}^{(0)}$
\begin{equation}\label{eq21}
A_{j} \Delta_{j}^{2} + 2B_{j} \Delta_{j} + C_{j} = O(\Delta_{j}^{3}),
\end{equation}
where
\begin{eqnarray}
\nonumber
A_{j} = \frac{1}{2}\left[\varphi''_{k}(c_{j}^{(0)})-\chi''_{j,0}(c_{j}^{(0)}) +
\tan\left(\chi_{j,0}(c_{j}^{(0)})\right)
\left(\varphi'_{k}(c_{j}^{(0)})-\chi'_{j,0}(c_{j}^{(0)})\right)^{2}\right],\\
\nonumber
B_{j} = \frac{1}{2}\left[\varphi'_{k}(c_{j}^{(0)})-\chi'_{j,0}(c_{j}^{(0)})\right], \quad
C_{j} = - \tan\left(\chi_{j,0}(c_{j}^{(0)})\right).
\end{eqnarray}
Hence, if the following conditions are fulfilled 
\begin{equation}\label{eq22}
\vert A_{j}C_{j}/B_{j}^{2}\vert < C_{\Delta}\ll 1,\quad
\vert B_{j}/A_{j}\vert > \delta_{0}/(1-C_{\Delta}), \quad
\vert C_{j}/2B_{j}\vert < \delta_{0}(1+C_{\Delta}).
\end{equation}
for some positive constant $C_{\Delta}$, the equation (\ref{eq21}) possesses a unique solution on the interval $(-\delta_{0}, \delta_{0})$.
\begin{lemmas}
For any $x\in I_{j,0}$ the derivatives of the function $\chi_{j,0}$ defined by (\ref{eq18}) admit the following estimates
\begin{equation*}
\chi'_{j,0}(x) = 
O({\rm  e}^{-2(1-\kappa_{0})\lambda_{0}/\varepsilon_{p,0}}),\quad
\chi''_{j,0}(x) =
O(\varepsilon^{-1}{\rm e}^{-(2-3\kappa_{0})\lambda_{0}/\varepsilon_{p,0}}).
\end{equation*}
\end{lemmas}
PROOF:  Let $\{\hat\varphi_{k}\}, \{\hat\lambda_{k}\}_{k=1}^{\infty}$ be two sequences of $C^{2}$-functions of a variable $x\in (a, b)$ such that $\hat\varphi''_{k}$, $\hat\lambda''_{k}$ are uniformly bounded and for any $x$ the sequences  
$\{\hat\varphi_{k}(x)\}, \{\hat\lambda_{k}(x)\}_{k=1}^{\infty}$ satisfy the conditions of 
Lemma 1. Then by Lemma 1 one may represent 
\begin{equation*}
A^{k}(x) = R(\theta_{k}(x)+\chi_{k}(x))Z(\mu_{k}(x))R(-\chi_{k}(x)).
\end{equation*}
Differentiating (\ref{eq8}) and using (\ref{eq9}) yields
\begin{equation}\label{eq23}
\chi'_{k}(x) = O\left(\delta^{-2} {\rm e}^{-2\lambda_{0}}\right),\quad
\chi''_{k}(x) = O\left(\delta^{-3} {\rm e}^{-2\lambda_{0}}\right).
\end{equation}
We apply (\ref{eq23}) to $M(A(x,\varepsilon),\tau_{0}-1)$ for $x\in I_{j,0}$. Then taking into account the specific dependence of functions $\varphi_{k}, \lambda_{k}$ on the parameter $\varepsilon$ (see (\ref{eq3})) one notes that  $\varphi'_{k}(x) = O(\varepsilon^{-1})$, $\varphi''_{k}(x) = O(\varepsilon^{-2})$. Together with the definition of $\delta_{0}$ this finishes the proof. $\square$

The non-degeneracy condition on the critical points of the function $\varphi_{k}$ and Lemma 3 imply the existence of positive constants $C_{1}, C_{2}$ such that
\begin{equation}\label{eq24}
\left\vert \frac{A_{j}C_{j}}{B_{j}^{2}}\right\vert \le 
C_{1}\frac{\varepsilon^{-1}{\rm e}^{-(2-\kappa_{0})\lambda_{0}/\varepsilon}}
{\left\vert \varphi'_{k}(c_{j}^{(0)})-\chi'_{j,0}(c_{j}^{(0)})\right\vert},\quad
\left\vert \frac{B_{j}}{4A_{j}\delta_{0}}\right\vert \ge 
C_{2}{\rm e}^{\kappa_{0}\lambda_{0}/\varepsilon}
\left\vert \varphi'_{k}(c_{j}^{(0)})-\chi'_{j,0}(c_{j}^{(0)})\right\vert.
\end{equation}
Hence, if 
\begin{equation}\label{eq25}
\left\vert \varphi'_{k}(c_{j}^{(0)})\right\vert \ge 
{\rm e}^{-\lambda_{0}\kappa_{0}\gamma/(1+\gamma)\varepsilon_{p,0}}
\end{equation}
we apply the inequalities (\ref{eq7}), (\ref{eq24}) and (\ref{eq13}) to conclude that the condition (\ref{eq22}) are fulfilled and the equation (\ref{eq20}) possesses a unique solution on the interval $I_{j,0}$. We denote by $\mathcal{E}''_{p,0}$ the set of those values of the parameter $\varepsilon$ for which the condition (\ref{eq25}) is violated. This set is a small neighborhood of one or both end-points of the interval $(\varepsilon_{p,1}, \varepsilon_{p,0}]$ and due to non-degeneracy of critical points of the functions $\varphi_{k}$ its Lebesgue measure is estimated as
\begin{equation*}
leb\left(\mathcal{E}''_{p,0}\right) = 
O\left({\rm e}^{-\lambda_{0}\kappa_{0}\gamma/(1+\gamma)\varepsilon_{p,0}}\right).
\end{equation*}
Introduce $\mathcal{E}_{p,1}=\hat{\mathcal{E}}_{p,1}\setminus \mathcal{E}''_{p,0}$. Then for any $\varepsilon\in \mathcal{E}_{p,1}$ there exists a unique solution of the equation (\ref{eq20}). Moreover, this solution denoted by $c_{j}^{(1)}(\varepsilon)$ satisfies the following estimate
$$
dist\left(c_{j}^{(1)}, c_{j}^{(0)}\right) =  O\left(\varepsilon{\rm e}^{-\left(2-\frac{1+2\gamma}{1+\gamma}\kappa_{0}\right)\lambda_{0}/\varepsilon_{p,0}}\right).
$$ 
These points constitute a set 
$$
\mathcal{C}_{1}(\varepsilon) = \bigcup\limits_{j=1}^{N(\varepsilon)}\{c_{j}^{(1)}(\varepsilon)\},
$$ 
which we call the first approximation of the critical set and define the layer of the order $1$ as
\begin{equation*}
L\mathcal{C}_{1} = 
\bigcup\limits_{j=1}^{N(\varepsilon)} (c_{j}^{(1)}-\delta_{1}, c_{j}^{(1)}+\delta_{1}) = \bigcup\limits_{j=1}^{N(\varepsilon)} I_{j,1}, \quad 
\delta_{1} = {\rm e}^{-\lambda_{0}\kappa_{1}/\varepsilon_{p,0}},
\end{equation*}
with some $\kappa_{1}$ which will be fixed later. Then in the same manner as for $L\mathcal{C}_{0}$ we can define collisions of the first order, their times $\tau_{j,j'}^{(1)}$ and
$$
\mathcal{E}'_{p,1}=\{\varepsilon\in \mathcal{E}_{p,1}: \exists\; (j,j'):\; \tau_{j,j'}^{(1)}<\tau_{1}\},
$$ 
where $\tau_{1}$ stands for the time of primary collisions of the first order.

Take the return time $q_{n_{J_{1}}}$ from the definition of condition $(A)$ and
define $\kappa_{1}$ as
$$
\frac{1}{q_{k_{J_{1}}+1}} < \delta_{1} = {\rm e}^{-\lambda_{0}\kappa_{1}/\varepsilon} = \frac{1}{C_{\omega}q_{k_{J_{1}}}^{1+\gamma}}.
$$
As at the previous step this choice implies 
\begin{eqnarray}
\nonumber
dist\left(\sigma_{\omega}^{m}(x), x\right) > \frac{1}{2 q_{n_{J_{1}}}}\gg 
\frac{1}{C_{\omega}q_{n_{J_{1}}}^{1+\gamma}}=\delta_{1},\; 0\le m\le q_{n_{J_{1}}}-1,\\
\nonumber
dist\left(\sigma_{\omega}^{q_{n_{J_{1}}}}(x), x\right) < \frac{1}{q_{n_{J_{1}}+1}} < 
\frac{1}{C_{\omega}q_{n_{J_{1}}}^{1+\gamma}}=\delta_{1}
\end{eqnarray}
and, hence, $\tau_{1} = q_{k_{J_{1}}}$. Analog of Lemma 2 gives
$$
leb(\mathcal{E}'_{p,1}) = 
O\left({\rm e}^{-\lambda_{0}\kappa_{1}\gamma/(1+\gamma)\varepsilon_{p,0}}\right).
$$

We exclude $\mathcal{E}'_{p,1}$ to get $\mathcal{E}_{p,2}=\mathcal{E}_{p,1}\setminus \mathcal{E}'_{p,1}$.
One may note that if $x\in I_{j,0}\setminus I_{j,1}$ then
\begin{equation*}
\Vert M(x, \tau_{0})\Vert \ge 
\left(C_{A}{\rm e}^{(1-\kappa_{0})\lambda_{0}/\varepsilon}\right)^{\tau_{0}-1}
\left(C_{A}{\rm e}^{(1-\kappa_{1})\lambda_{0}/\varepsilon}\right)=
C_{A}^{\tau_{0}}{\rm e}^{\left(\tau_{0}-(\tau_{0}-1)\kappa_{0}-\kappa_{1}\right)\lambda_{0}/\varepsilon}.
\end{equation*}
Let $\varepsilon\in \mathcal{E}_{p,2}$ and $x\in I_{j,1}$, then one notes that for $1\le m\le \tau_{1}-1$ the finite trajectory $\{\sigma_{\omega}^{m}\}$ cannot fall into $I_{j,1}$, but only into $I_{j,0}\setminus I_{j,1}$. Taking this into account we apply successively (\ref{eq8}),  (\ref{eq9}), (\ref{eq11}) together with the first inequality (\ref{eq12}) to represent $M(A(x,\varepsilon),\tau_{1}-1)$ as
\begin{equation*}
M(A(x,\varepsilon),\tau_{1}-1) = 
R(\theta_{j,1}(x)+\chi_{j,1}(x))Z(\mu_{j,1}(x))R(-\chi_{j,1}(x)),
\end{equation*}
where 
\begin{multline}\label{eq26}
\Bigl\vert\chi_{j,1}(x) - \chi_{j,0}(x)\Bigr\vert \le \varepsilon 
{\rm e}^{-2(\tau_{0}-1)(1-\kappa_{0})\lambda_{0}/\varepsilon + \log(2 q_{n_{J_{0}}+1})},\\
\mu_{j,1}(x)\ge \frac{1}{\varepsilon}\left((\tau_{1}-1)(\lambda_{0}-\kappa_{0}-\log C_{A}) - \sum\limits_{k=n_{J_{0}}}^{n_{J_{1}}-1}\log(2 q_{k+1})\right).
\end{multline}
Finally, due to the Brjuno's condition one gets
\begin{eqnarray}
\nonumber
\Bigl\vert\chi_{j,1}(x) - \chi_{j,0}(x)\Bigr\vert = O(\varepsilon{\rm e}^{-2(\tau_{0}-1)\left[(1-\kappa_{0})\lambda_{0}/\varepsilon - C_{B}/2\right]}),\\
\nonumber
\mu_{j,1}(x)\ge \frac{\tau_{1}-1}{\varepsilon}\Bigl(\lambda_{0}-\kappa_{0}-\log C_{A} - C_{B}\Bigr).
\end{eqnarray}
Similarly to (\ref{eq20}) we consider an equation
\begin{equation}
\cos\left(\varphi_{k}(x) - \chi_{j,1}(x)\right) = 0,\quad x\in I_{j,1},
\end{equation}
with the same $k$ defined by (\ref{eq19}). Since (\ref{eq25}) holds for 
$\varepsilon\in \mathcal{E}_{p,2}$ and due to the inequality (\ref{eq7}) this equation has a unique solution in a small vicinity of $c_{j}^{(1)}$. We denote  this solution by $c_{j}^{(2)}(\varepsilon)$. Then it satisfies the following estimate
\begin{equation*}
dist\left(c_{j}^{(2)}, c_{j}^{(1)}\right) =  O\left(\varepsilon
{\rm e}^{-2(\tau_{0}-1)\left[(1-\kappa_{0})\lambda_{0}/\varepsilon - C_{B}\right] + 
\lambda_{0}\kappa_{0}\gamma/(1+\gamma)\varepsilon}\right).
\end{equation*}
We proceed in the same manner. Assuming the $n$-th approximation of the critical set
$$
\mathcal{C}_{n}(\varepsilon) = \bigcup\limits_{j=1}^{N(\varepsilon)}\{c_{j}^{(n)}(\varepsilon)\},
$$ 
is constructed one may define the layer of order $n$ 
\begin{equation*}
L\mathcal{C}_{n} = 
\bigcup\limits_{j=1}^{N(\varepsilon)} (c_{j}^{(n)}-\delta_{n}, c_{j}^{(n)}+\delta_{n}) = \bigcup\limits_{j=1}^{N(\varepsilon)} I_{j,n}, \quad 
\delta_{n} = {\rm e}^{-\lambda_{0}\kappa_{n}/\varepsilon_{p,0}},
\end{equation*}
where $\delta_{n}$ is such that
$$
\frac{1}{q_{k_{J_{n}}+1}} < \delta_{n} = {\rm e}^{-\lambda_{0}\kappa_{n}/\varepsilon} = \frac{1}{C_{\omega}q_{k_{J_{n}}}^{1+\gamma}}.
$$
Dynamics of the layer $L\mathcal{C}_{n}$ under the rotation $\sigma_{\omega}$ gives rise to the notion of collisions of the $n$-th order and their times $\tau_{j,j'}^{(n)}$. We also construct the set 
$$
\mathcal{E}'_{p,n}=\{\varepsilon\in \mathcal{E}_{p,n}: \exists\; (j,j'):\; \tau_{j,j'}^{(n)}<\tau_{n}\},
$$ 
where $\tau_{n}$ stands for the time of primary collisions of the $n$-th order. A corresponding modification of 
Lemma 2 gives
\begin{equation}\label{eq27}
leb(\mathcal{E}'_{p,n}) = 
O\left({\rm e}^{-\lambda_{0}\kappa_{n}\gamma/(1+\gamma)\varepsilon_{p,0}}\right).
\end{equation} 
Eliminating $\mathcal{E}'_{p,n}$ we obtain the $n+1$-th approximation of a set of "good" values of the parameter $\varepsilon$, namely, $\mathcal{E}_{p,n+1}=\mathcal{E}_{p,n}\setminus \mathcal{E}'_{p,n}$. For $\varepsilon\in \mathcal{E}_{p,n+1}$ and $x\in I_{j,n}$ by similar arguments used to get (\ref{eq26}) one may represent $M(A(x,\varepsilon),\tau_{n}-1)$ in the following form
\begin{equation*}
M(A(x,\varepsilon),\tau_{n}-1) = 
R(\theta_{j,n}(x)+\chi_{j,n}(x))Z(\mu_{j,n}(x))R(-\chi_{j,n}(x)),
\end{equation*}
where 
\begin{multline*}
\Bigl\vert\chi_{j,n}(x) - \chi_{j,n-1}(x)\Bigr\vert \le \varepsilon 
{\rm e}^{-2(\tau_{n-1}-\tau_{n-2}-1)(1-\kappa_{0})\lambda_{0}/\varepsilon} \cdot
\prod\limits_{k=n_{J_{n-2}}+1}^{n_{J_{n-1}}}2 q_{k+1},\\
\mu_{j,n}(x)\ge \frac{1}{\varepsilon}\left((\tau_{n}-1)(\lambda_{0}-\kappa_{0}-\log C_{A}) - \sum\limits_{k=n_{J_{n-1}}}^{n_{J_{n}}-1}\log(2 q_{k+1})\right).
\end{multline*}
Using the Brjuno's condition we obtain
\begin{eqnarray}\label{eq28}
\nonumber
\Bigl\vert\chi_{j,n}(x) - \chi_{j,n-1}(x)\Bigr\vert = O(\varepsilon
{\rm e}^{-2(\tau_{n-1}-\tau_{n-2}-1)\left[(1-\kappa_{0})\lambda_{0}/\varepsilon - C_{B}/2\right]}),\\
\mu_{j,n}(x)\ge \frac{\tau_{n}-1}{\varepsilon}\Bigl(\lambda_{0}-\kappa_{0}-\log C_{A} - C_{B}\Bigr).
\end{eqnarray}
This leads to definition of $c_{j}^{(n+1)}$ as a unique solution of the equation
\begin{equation*}
\cos\left(\varphi_{k}(x) - \chi_{j,n}(x)\right)=0,\quad x\in I_{j,n}
\end{equation*}
which satisfies the bound
\begin{equation}\label{eq29}
dist\Bigl( c_{j}^{(n+1)}, c_{j}^{(n)}\Bigr) = O\left(\varepsilon
{\rm e}^{-2(\tau_{n-1}-\tau_{n-2}-1)\left[(1-\kappa_{0})\lambda_{0}/\varepsilon - C_{B}\right]  + \lambda_{0}\kappa_{0}\gamma/(1+\gamma)\varepsilon}\right).
\end{equation} 
The existence and uniqueness of this solution follows from the inequalities (\ref{eq7}), (\ref{eq25}) and (\ref{eq28}).

Finally, passing to the limit $n\to \infty$ one may note that the estimates (\ref{eq27}), (\ref{eq29}) yield the following lemma
\begin{lemmas}
There exist limit sets 
$$
\mathcal{C}_{\infty}(\varepsilon) = \bigcup\limits_{j=1}^{N(\varepsilon)}\{c_{j}^{(\infty)}(\varepsilon)\},\quad
c_{j}^{(\infty)} = \lim\limits_{n\to \infty} c_{j}^{(n)}
$$ 
and
$$
\mathcal{E}_{p,\infty} = \mathcal{E}_{p,0}\setminus 
\bigcup\limits_{n=0}^{\infty} \mathcal{E}'_{p,n}
$$
such that the Lebesgue measure of $\mathcal{E}_{p,0}\setminus \mathcal{E}_{p,\infty}$ admits an estimate
$$ 
leb(\mathcal{E}_{p,0}\setminus \mathcal{E}_{p,\infty}) = 
O\left({\rm e}^{-\lambda_{0}\kappa_{0}\gamma/(1+\gamma)\varepsilon_{p,0}}\right).
$$
\end{lemmas}

\section{Lower bound for the Lyapunov exponent}
\setcounter{equation}{0}

In this subsection we always suppose that $\varepsilon\in \mathcal{E}_{p,\infty}$.
\begin{defs}
We say that $x\in \mathbb{T}^{1}$ has property $(H)$ if
$$
\begin{cases}
\sigma_{\omega}^{m}(x)\notin L\mathcal{C}_{0},\quad \forall\; 0\le m < \tau_{0},\\
\sigma_{\omega}^{m}(x)\notin L\mathcal{C}_{n},\quad \forall\; \tau_{n-1}\le m < \tau_{n}.
\end{cases}
$$
\end{defs}
Denote the set of all $x\in \mathbb{T}^{1}$ which possess the property $(H)$ by $X_{h}$ and its complement by $X_{e} = \mathbb{T}^{1}\setminus X_{h}$. Then definition of $\delta_{n}$ together with the property $(A)$ imply 
\begin{multline*}
leb\left(X_{e}\right)\le 2p\sum\limits_{n=0}^{\infty}\tau_{n} \delta_{n}\le 
2p\sum\limits_{n=0}^{\infty}C_{\omega}^{-\frac{1}{1+\gamma}} 
\delta_{n}^{\frac{\gamma}{1+\gamma}}\le \\
2pC_{\omega}^{-\frac{1}{1+\gamma}}\sum\limits_{n=0}^{\infty} 
\left(C_{\delta}^{\frac{\gamma}{1+\gamma}}\right)^{n}
\delta_{0}^{\frac{\gamma}{1+\gamma}}=
\frac{2p C_{\omega}^{-\frac{1}{1+\gamma}}}{1-C_{\delta}^{\frac{\gamma}{1+\gamma}}}
{\rm e}^{-\lambda_{0}\kappa_{0}\gamma/(1+\gamma)\varepsilon_{p,0}}.
\end{multline*}
Taking into account that $p=O(\varepsilon_{p,0}^{-1})$ one concludes
\begin{equation*}
leb(X_{h}) = 1 - O\left(\varepsilon_{p,0}^{-1}
{\rm e}^{-\lambda_{0}\kappa_{0}\gamma/(1+\gamma)\varepsilon_{p,0}}\right).
\end{equation*}
Hence, the Lebesgue measure of $X_{h}$ is positive. Since the  Lyapunov exponent exists and
is constant a.e., it is sufficient to estimate $\Lambda(x)$ for $x\in X_{h}$.

Consider $x\in X_{h}$ and its trajectory $x_{m} = \sigma_{\omega}^{m}(x), m=0, 1, \dots$. We denote by $m_{i}$ the first time when 
$x_{m_{i}}\in L\mathcal{C}_{n_{i}}\setminus L\mathcal{C}_{n_{i}+1}$. By construction,
$m_{i+1} - m_{i}\ge \tau_{i}$. Let $m_{k}\le m < m_{k+1}$. Then, using the property $(H)$ and the Brjuno's condition, one gets
\begin{multline}
\nonumber
\frac{1}{m}\log\Vert M(x, m)\Vert \ge 
\frac{1}{m}\log\left(C_{A}{\rm e}^{(1-\kappa_{0})\lambda_{0}/\varepsilon}\right)^{m}\cdot
\prod\limits_{i=1}^{m_{k}}\delta_{n_{i}+1}\ge\\
(1 - \kappa_{0})\lambda_{0}/\varepsilon +\log C_{A}-\frac{1}{m}\sum\limits_{i=1}^{m_{k}}\log 2q_{n_{i}+1}\ge
(1 - \kappa_{0})\lambda_{0}/\varepsilon +\log C_{A}-C_{B}.
\end{multline}
Thus, we arrive at the following lemma
\begin{lemmas}
If $\varepsilon\in \mathcal{E}_{p,\infty}$ the integrated Lyapunov exponent
\begin{equation*}
\Lambda_{0}\ge (1 - \kappa_{0})\lambda_{0}/\varepsilon +\log C_{A}-C_{B}.
\end{equation*}
\end{lemmas}

Summation over all $p$ proves the first and the second assertion of Theorem 4. However, one may note that if one takes a point $x\in \mathbb{T}^{1}$ sufficiently close to the critical set $\mathcal{C}_{\infty}(\varepsilon)$ the Lyapunov exponent $\Lambda(x)$ can be made arbitrarily small. Thus, in the case of non-constant functions $\varphi_{k}$ the cocycle cannot be reducible for suficiently large set of the parameters and, hence, does not possess the exponential dichotomy.

{\bf Remark} We emphasize that the Brjuno's condition and the condition $(A)$ play different roles in the paper. While the Brjuno's condition is mainly used to get the lower bound for the Lyapunov exponent, the condition $(A)$ is applied to avoid the secondary collisions and to estimate the Lebesgue measure of the set $\mathcal{E}'_{p,n}$. Omitting the condition $(A)$ (or a condition similar to it) leads to growth of $leb(\mathcal{E}'_{p,n})$. In this case the measure of $\mathcal{E}'_{p,n}$ will not be exponentially small, however, it tends to zero with decreasing of $\varepsilon$ \cite{LSY}.

\section*{Acknowledgements}

The research was supported by RFBR grant (project No. 17-01-00668/19).


\end {document}